\documentclass[11pt]{article}
\usepackage{latexsym}
\usepackage{amsmath, amsthm}
\usepackage[all,dvips,ps]{xy}
\usepackage{amsfonts}
\usepackage{mathrsfs}

\newtheorem{lem}{Lemma}[section]
\newtheorem{thm}{Theorem}[section]

\newcommand{\Diag}{\mathrm{Diag}\,}

\newcommand{\Rs}{\mathbb{R}}

\newcommand{\bz}{{\bf 0} }

\newcommand{\bpr}{{\bf Proof.} \hspace{1 em}}
\newcommand{\epr}{ \\ \hspace*{4.5in} $\Box$ }
\newcommand{\beq}{ \begin{equation} }
\newcommand{\eeq}{ \end{equation} }
\newcommand{\bt}{ \begin{tabular} }
\newcommand{\et}{ \end{tabular} }

\begin{document}

\bibliographystyle{plain}
\title{Graph Connectivity and Universal Rigidity of Bar Frameworks
 \thanks{Research supported by the Natural Sciences and Engineering
         Research Council of Canada.} }
\vspace{0.3in}
        \author{ A. Y. Alfakih
  \thanks{E-mail: alfakih@uwindsor.ca}
  \\
          Department of Mathematics and Statistics \\
          University of Windsor \\
          Windsor, Ontario N9B 3P4 \\
          Canada
}

\date{Jul 9, 2014. Revised \today}
\maketitle


\vspace{0.1in}

\begin{abstract}
Let $G$ be a graph on $n$ nodes.
In this note, we prove that if $G$ is $(r+1)$-vertex connected, $1 \leq r \leq n-2$, 
then there exists a configuration $p$
in general position in $\Rs^r$ such that the bar framework $(G,p)$ is universally rigid.
The proof is constructive, and is based on a theorem by Lov\'{a}sz {\em et al}
concerning orthogonal representations and connectivity of graphs \cite{lss89, lss00}.
\end{abstract}

\section{Introduction}
Let $G=(V,E)$ be a graph on $n$ nodes. 
$G$ is said to be {\em $k$-vertex connected}, or simply {\em $k$-connected},
if $n=k+1$ and $G$ is the complete graph, or if $n \geq k+2$ and 
there does not exist a set of $(k-1)$ nodes whose deletion disconnects $G$.
A {\em bar framework} in $\Rs^r$
is a simple incomplete connected graph $G$ whose nodes are points
$p^1,\ldots,p^n$ in $\Rs^r$; and whose edges are line segments, each joining a pair
of these points.
The points $p^1, \ldots, p^n$ will be denoted collectively by $p$, and
the bar framework will be denoted by $(G,p)$. Also, we will refer to $p$ as
the {\em configuration} of the bar framework.
A configuration $p$ (or a framework $(G,p)$) is {\em $r$-dimensional} if the points
$p^1,\ldots,p^n$ affinely span $\Rs^r$.
Moreover, a configuration $p$ (or a framework $(G,p)$) is in general position in
$\Rs^r$ if every $r+1$ points in configuration $p$ are affinely independent.

An $r'$-dimensional bar framework  $(G,p')$ is {\em equivalent} to an $r$-dimensional bar
framework $(G,p)$ if:
\beq \label{defd}
||{p'}^i-{p'}^j||^2 = ||p^i - p^j ||^2 \quad \text{for each $\{i,j\} \in E(G)$},
\eeq
where $||.||$ denotes the Euclidean norm and $E(G)$ denotes the edge set of graph $G$.
On the other hand, two $r$-dimensional bar frameworks $(G,p)$ and $(G,p')$
are {\em congruent} if:
\beq \label{defd}
||{p'}^i-{p'}^j||^2 = ||p^i - p^j ||^2 \quad \text{for all $i,j=1,\ldots,n$}.
\eeq
An $r$-dimensional bar framework $(G,p)$ is said to be {\em universally rigid} if
there does not exist an $r'$-dimensional bar framework $(G,p')$, where $r'$ is a positive integer
$\leq n-1$, such that $(G,p')$ is equivalent but not congruent to $(G,p)$.

An immediate necessary condition for an $r$-dimensional bar framework $(G,p)$ on $n$ nodes $(r \leq n-2$)
in general position in $\Rs^r$ to be universally rigid is
that graph $G$ should be $(r+1)$-connected \cite{hen92}.  For suppose that $G$ is not $(r+1)$-connected. Then
there exists a set of $r$ nodes, say $X$, whose removal disconnects $G$.
Let $V(G)= V_1 \cup X \cup V_2$ be a partition of the nodes of $G$, where $V_1$ and $V_2$ are non-empty,
such that there are no edges joining nodes in $V_1$ to nodes in $V_2$.
The points $\{p^i: i \in X\}$ lie in a hyperplane $H$ in $\Rs^r$, and the points 
$\{p^i: i \in  V_1 \cup V_2\}$ do not lie in $H$ since $p$ is in general position in $\Rs^r$.
For all nodes $i \in V_2$, let $q^i$
be the reflection of $p^i$ with respect to $H$ and let
$p'=\{ p^i: i \in (V_1 \cup X)\} \cup \{q^i: i \in V_2\}$. Thus $(G,p')$ is an $r$-dimensional  bar
framework that is equivalent but not congruent to $(G,p)$, and hence $(G,p)$ is not universally rigid.
This raises the question of whether the assumption of $(r+1)$-connectivity of graph $G$ alone is sufficient for
the existence of some $r$-dimensional configuration $p$ in general position in $\Rs^r$ such that
the bar framework $(G,p)$ is universally rigid.
The following theorem, which is our main result, is an affirmative answer to this question.

\begin{thm} \label{thmmain}
Let $G$ be a graph on $n$ nodes and assume that $G$ is $(r+1)$-vertex connected, where $1 \leq r \leq n-2$.
Then
there exists an $r$-dimensional bar framework $(G,p)$
in general position in $\Rs^r$ such that $(G,p)$ is universally rigid.
\end{thm}
The proof of Theorem \ref{thmmain}, which is given in Section 3,
 is constructive and is based on a theorem by Lov\'{a}sz {\em et al}
\cite{lss89,lss00} concerning orthogonal representations and connectivity of graphs.

Note that the complete bipartite graph $K_{3,3}$ is $3$-connected. Thus,
Theorem \ref{thmmain} provides a negative answer to a question raised by Yinyu Ye as to whether
in every universally rigid $2$-dimensional bar framework $(G,p)$ in general position in $\Rs^2$,
graph $G$ must contain a triangle.

\section{Preliminaries}

This section presents the necessary mathematical background. The first subsection reviews basic
definitions and results on stress and Gale matrices, and their role in the problem of universal rigidity.
The second subsection focuses on vertex connectivity and orthogonal representations of graphs.

\subsection{Stress and Gale Matrices}

Stress matrices play a key role in the study of universal rigidity.
An {\em equilibrium stress} (or simply a {\em stress}) of a bar framework $(G,p)$ is a real-valued function
$\omega$ on $E(G)$ such that:
\beq \label{defw}
\sum_{j:\{i,j\} \in E(G)} \omega_{ij} (p^i - p^j) = \bz \mbox{ for each } i=1,\ldots,n.
\eeq
Here we use the bold zero ``$\bz$" to denote the zero vector or the zero matrix of appropriate dimensions.
Let $E(\overline{G})$ denote the edge set of graph $\overline{G}$, the complement graph of $G$. i.e.,
\[
E(\overline{G})= \{ \{i,j\}: i \neq j , \{i,j\} \not \in E(G) \},
\]
and let $\omega =(\omega_{ij})$ be a stress of $(G,p)$. Then the
$n \times n$ symmetric matrix $\Omega$ where
\beq \label{defO}
\Omega_{ij} = \left\{ \begin{array}{ll} -\omega_{ij} & \mbox{if } \{i,j\} \in E(G), \\
                        0   & \mbox{if }  \{i,j\}  \in E(\overline{G}), \\
                   {\displaystyle \sum_{k:\{i,k\} \in E(G)} \omega_{ik}} & \mbox{if } i=j,
                   \end{array} \right.
\eeq
is called the {\em stress matrix
associated with $\omega$}, or a {\em stress matrix} of $(G,p)$. Sufficient and necessary conditions, in terms
of stress matrices, for universal rigidity of bar frameworks are discussed in 
\cite{alf07a, con99, con82, alf10, gt14}.
The first sufficient condition for universal rigidity under the assumption that configuration
$p$ is in general position was given in \cite{ay13}.

\begin{thm}[Alfakih and Ye \cite{ay13}]
\label{thmay}
Let $(G,p)$ be an $r$-dimensional bar framework on $n$ nodes in $\Rs^r$, for some $r
\leq n-2$. If the
following two conditions hold:
\begin{enumerate}
\item There exists a positive semidefinite stress matrix $\Omega$ of $(G,p)$ of rank
$n-r-1$,
\item The configuration $p$ is in general position.
\end{enumerate}
Then $(G,p)$ is universally rigid.
\end{thm}

Theorem \ref{thmay} was generalized and strengthened in \cite{an13}, but it will suffice for the purposes
of this note.

Stress matrices are intimately related to Gale matrices and Gale transform \cite{gal56,gru67}.
This relation is a crucial step in connecting stress matrices to orthogonal representations
of graphs.
Given an $r$-dimensional bar framework $(G,p)$ in $\Rs^r$, let
\beq \label{defP}
P:=
\left[ \begin{array}{c} (p^1)^T  \\ \vdots  \\ (p^n)^T \end{array} \right].
\eeq
Then $P$ is called the {\em configuration matrix} of $p$ (or of framework $(G,p)$). Moreover, let
\beq \label{defz}
Z:=
\left[ \begin{array}{c} (z^1)^T  \\ \vdots  \\ (z^n)^T \end{array} \right]
\eeq
be any $n \times (n-r-1)$ matrix
whose columns form a basis of the null space of the matrix
\beq \label{defPe}
\left[ \begin{array}{cc} P^T \\ e^T \end{array} \right],
\eeq
where $e$ is the vector of all 1's in $\Rs^n$.
Note that the matrix in (\ref{defPe}) has full row rank since $(G,p)$ is $r$-dimensional.
Then
$Z$ is called a {\em Gale matrix} of configuration $p$ (or of bar framework $(G,p)$), and
$z^1,\ldots,z^n$ are called, respectively, {\em Gale transforms} of $p^1,\ldots,p^n$.
Note that $z^1,\ldots,z^n$ are vectors in $\Rs^{n-r-1}$. Also, note that if $Z$ is a Gale
matrix of $(G,p)$ and $Q$ is any nonsingular matrix of order $n-r-1$, then $ZQ$ is also
a Gale matrix of $(G,p)$.

It readily follows from (\ref{defO}) that a stress matrix $\Omega$ satisfies the following
equations:
\[
\Omega P = \bz, \Omega \, e =\bz \mbox{ and $\Omega_{ij}=0$ for each } \{i,j\} \in E(\overline{G}) .
\]
Thus, the columns of $\Omega$ belong to the null space of the matrix in (\ref{defPe}).
Accordingly, we have  the following lemma.

\begin{lem}[Alfakih \cite{alf10}] \label{lemOZ}
Let $(G,p)$ be an $r$-dimensional bar framework on $n$ nodes in $\Rs^r$, $r \leq n-2$, and let
$Z$ be a Gale matrix of $(G,p)$. Further, let $\Omega =  Z \Psi Z^T$ for some symmetric 
matrix $\Psi$ of order $n-r-1$. If
\beq
(Z \Psi Z^T)_{ij}= 0 \mbox{ for each } \{i,j\} \in E(\overline{G}),
\eeq
then $\Omega$ is a stress matrix of $(G,p)$ and
rank $\Omega$ = rank $\Psi \leq n-r-1$.
\end{lem}
A point worth observing here is that, by Lemma \ref{lemOZ}, $\Omega$ is positive semidefinite with rank $n-r-1$
if and only if $\Psi$ is positive definite if and only if
there exist Gale transforms ${z}'^1,\ldots,{z'}^n$ of $p^1,\ldots,p^n$ such that $({z'}^i)^T {z'}^j=0$ for each
$\{i,j\} \in E(\overline{G})$. 

The usefulness of Gale transform in the study of universal rigidity under the general position assumption stems from
the fact that Gale matrix $Z$, and consequently Gale transforms $z^1,\ldots,z^n$, encode the affine dependence
of points $p^1, \ldots, p^n$.

\begin{lem}
\label{lemgp}
Let $(G,p)$ be an $r$-dimensional bar framework on $n$ nodes in $\Rs^r$ and let $z^i$ be a Gale transform of
$p^i$ for $i=1,\ldots,n$. Then $(G,p)$ is in general position if and only if
any size-$(n-r-1)$ subset of $\{z^1,\ldots,z^n\}$ is linearly independent; i.e.,
any $(n-r-1) \times (n-r-1)$ submatrix of Gale matrix $Z$ is nonsingular.
\end{lem}
For a proof of Lemma \ref{lemgp} see e.g. \cite{alf07a}.

\subsection{Graph Connectivity and Orthogonal Representations}

An orthogonal representation of a graph $G$ in $\Rs^k$
is a mapping of each node $i$ of $G$ into a vector $x^i$ in
$\Rs^k$ such that  $x^i$ is orthogonal to $x^j$ for every pair of nonadjacent nodes $i$ and $j$ of $G$;
i.e., $(x^i)^Tx^j=0$ for each $\{i,j\} \in E(\overline{G})$.
The vectors $x^1,\ldots,x^n$ are called the {\em representing vectors}.
Orthogonal representations of graphs were introduced by Lov\'{a}sz in \cite{lov79} in his study of
the Shannon capacity of a graph. Obviously, $x^i=\bz$ for each node $i$ of $G$ is a trivial
orthogonal representation of $G$. Thus, in order to exclude such degenerate cases,
orthogonal representations are required to satisfy the condition that
any size-$k$ subset of $\{x^1,\ldots, x^n\}$ is linearly independent.
The following theorem by Lov\'{a}sz {\em et al} is crucial for this note.

\begin{thm}[Lov\'{a}sz {\em et al} \cite{lss89,lss00}] \label{thmlss}
Let $G$ be a graph on $n$ nodes, then $G$ is $(r+1)$-vertex connected, $r \leq n-2$, if and only if
$G$ has an orthogonal representation in $\Rs^{n-r-1}$ such that every size-($n-r-1$) subset
of the representing vectors is linearly independent.
\end{thm}

Let $G$ be a graph on $n$ nodes such that
each node of $G$ has a degree at least $r+1$; i.e., each node of $G$ has at most
$n-r-2$ non-adjacent nodes.
Lov\'{a}sz {\em et al} \cite{lss89,lss00} presented the following simple randomized
algorithm to construct an orthogonal representation of $G$ in $\Rs^{n-r-1}$.
Fix an ordering $(1,\ldots,n)$ of the nodes of $G$. Then the representing vectors
$x^1,\ldots,x^n$ are selected sequentially as follows. Select $x^1$ to be a uniformly random unit vector
in $\Rs^{n-r-1}$.  For $j=2,\ldots,n$, having selected $x^1,\ldots, x^{j-1}$, select $x^j$
to be a uniformly random unit vector from the subspace of $\Rs^{n-r-1}$ that is orthogonal
to the span of $\{x^i: i < j \mbox{ and } \{i,j\} \in E(\overline{G}) \}$. This subspace has dimension
$\geq 1$ since the dimension of the span of $\{x^i: \{i,j\} \in E(\overline{G})\}$ is
$\leq n-r-2$.
Now Lov\'{a}sz {\em et al} proved that if, in addition, $G$ is $(r+1)$-connected, then,
with probability 1, the orthogonal representation
constructed by this algorithm has the property that
every size-$(n-r-1)$ subset of $\{x^1,\ldots,x^n\}$ is linearly independent.

\section{Proof of Theorem \ref{thmmain}}

Let $G$ be a graph on $n$ nodes and assume that $G$ is $(r+1)$-connected.
Then by Theorem \ref{thmlss} there exist vectors $x^1,\ldots,x^n$ in $\Rs^{n-r-1}$ such that
$(x^i)^Tx^j=0$ for each $\{i,j\} \in E(\overline{G})$; and
every size-$(n-r-1)$ subset of $\{x^1,\ldots,x^n\}$ is linearly independent.

Let $X^T$ be the $(n-r-1) \times n$ matrix whose $i$th column is equal to $x^i$; i.e.,
\beq \label{defX}
X^T=\left[ \begin{array}{ccc} x^1  \; x^2 \; \cdots \; x^n \end{array} \right].
\eeq
 Then $(XX^T)_{ij} = (x^i)^T x^j = 0$ for  each $\{i,j\} \in E(\overline{G})$ and any square submatrix of $X^T$ of
 order $n-r-1$ is nonsingular.
The following two simple lemmas are needed.

\begin{lem} \label{lemxi}
Let $X^T$ be the matrix defined in (\ref{defX}), then there exists a vector  $\xi=(\xi_i)$ in $\Rs^n$
such that $X^T \xi = \bz$ and $\xi_i \neq 0$ for each $i=1,\ldots,n$.
\end{lem}

\bpr
Without loss of generality let $\left[ \begin{array}{c} I_{r+1} \\ B \end{array} \right]$ be
the $n \times (r+1)$ matrix
whose columns form a basis of the null space of $X^T$, where $I_{r+1}$ denotes the identity matrix
of order $r+1$, and $B$ is an $(n-r-1) \times (r+1)$ matrix.
Then each column of $\left[ \begin{array}{c} I_{r+1} \\ B \end{array} \right]$
has exactly $r$ zero entries; i.e., $B$ has no zero entries.
For suppose that $B$ has a zero entry, say $B_{1 1}=0$. Then
the size-$(n-r-1)$ set $\{x^1, x^{r+3}, x^{r+4}, \ldots, x^n\}$ is linearly dependent, a contradiction.
Therefore, $\xi$ is obtained be an appropriate linear combination of the columns of
$\left[ \begin{array}{c} I_{r+1} \\ B \end{array} \right]$.
\epr

\begin{lem} \label{lem2}
Let $X$ and $\xi$ be as in Lemma \ref{lemxi} and let $Z=\Diag(\xi) X$, where
$\Diag(\xi)$ is the diagonal matrix formed from the vector $\xi$. Furthermore,
Let $P$ be the $n \times r$ matrix whose columns form a basis of the null space
of
\beq \label{def2P}
\left[ \begin{array}{c} Z^T \\ e^T \end{array} \right],
\eeq
and let $p$ be the configuration in $\Rs^r$ whose configuration matrix is $P$.
Then $p$ is in general position in $\Rs^r$ and $Z$ is
a Gale matrix of $p$.
\end{lem}

\bpr
Note that $Z$ is $n \times (n-r-1)$ and $Z^Te=X^T \xi = \bz$.
Since $\xi$ has no zero entries, the matrix $\Diag(\xi)$ is nonsingular. Thus
$Z$ has full column rank and hence, by the definition of $P$ in the lemma,
it follows that $Z$ is a Gale matrix of configuration $p$. Furthermore,
every square submatrix of $Z$ of order $n-r-1$ is nonsingular. Therefore,
by Lemma \ref{lemgp}, configuration $p$ is in general position in $\Rs^r$.
\epr

To complete the proof of Theorem \ref{thmmain}, let $\Omega=ZZ^T= \Diag(\xi)XX^T \Diag(\xi)$.
Then, obviously, $\Omega$ is positive
semidefinite of rank $n-r-1$. Moreover, let $\{i,j\} \in E(\overline{G})$, then
$\Omega_{ij}=\xi_i \xi_j (XX^T)_{ij}=0$. Hence, $\Omega$ is a stress matrix of
the $r$-dimensional framework $(G,p)$ in $\Rs^r$ whose configuration matrix $P$ is as given in Lemma \ref{lem2}.
Since $(G,p)$ is in general position, it follows from Theorem \ref{thmay} that
$(G,p)$ is universally rigid. 

\noindent{\bf \large Acknowledgements} The author would like to thank Tibor Jord\'{a}n for useful comments on an earlier
version of this note.


\end{document}